\documentclass[preprint]{elsarticle}

\usepackage{subfig}
\usepackage{hyperref}
\usepackage{amsmath,amssymb,fixmath}
\usepackage{epsfig,epstopdf}
\usepackage{booktabs}
\usepackage{multirow}
\journal{Noname}
\begin{document}

\begin{frontmatter}

\title{Comments on \textquotedblleft Momentum fractional LMS for power signal parameter estimation\textquotedblright}


\author[mymainaddress]{Shujaat Khan}
\ead{shujaat@kaist.ac.kr}

\author[b,c]{Imran Naseem\corref{mycorrespondingauthor}}
\ead{imran.naseem@uwa.edu.au}

\author[c]{Alishba Sadiq}
\ead{alishba.sadiq@pafkiet.edu.pk}

\author[d]{Jawwad Ahmad}
\ead{jawwad@uit.edu}

\author[e,f]{Muhammad Moinuddin}
\cortext[mycorrespondingauthor]{Corresponding author}
\ead{mmsansari@kau.edu.sa}

\address[mymainaddress]{Department of Bio and Brain Engineering, Korea Advanced Institute of Science and Technology (KAIST), Daejeon, Republic of Korea.}
\address[b]{School of Electrical, Electronic and Computer Engineering, The University of Western Australia, 35 Stirling Highway, Crawley, Western Australia 6009, Australia.}
\address[c]{College of Engineering, Karachi Institute of Economics and Technology,	Korangi Creek, Karachi 75190, Pakistan.}
\address[d]{Department of Electrical Engineering, Usman Institute of Technology (UIT), Karachi, Pakistan.}
\address[e]{Center of Excellence in Intelligent Engineering Systems (CEIES), King Abdulaziz University, Jeddah, Saudi Arabia.}
\address[f]{Electrical and Computer Engineering Department, King Abdulaziz University, Jeddah, Saudi Arabia.}

\begin{abstract}
The purpose of this paper is to indicate that the recently proposed Momentum fractional least mean squares (mFLMS) algorithm has some serious flaws in its design and analysis.  Our apprehensions are based on the evidence we found in the derivation and analysis in the paper titled: \textquotedblleft \textit{Momentum fractional LMS for power signal parameter estimation}\textquotedblright.  In addition to the theoretical bases our claims are also verified through extensive simulation results.  The experiments clearly show that the new method does not have any advantage over the classical least mean square (LMS) method.
\end{abstract}


\end{frontmatter}


\section{Introduction}\label{Sec:Intro}
The least mean square (LMS) is one of the most widely used algorithms in adaptive signal processing \cite{LMSBook1}.  It has a number of variants to deal with various signals and environmental conditions \cite{ECGLMS, NCLMS, Douglas1994, CLMS, LMS11, QKLMS, ECLMS, TSPLMS, qLMS, FCLMS}.  To improve the convergence performance of the conventional LMS, different methods have been proposed based on adaptive step-size notion \cite {VSS, MnZ, RVSS1, sulyman2003convergence}.  Amongst all the variants, an important modification is the one proposed by J. G. Proakis and is called as the momentum LMS (mLMS) \cite{proakis1974channel}.  Unlike the conventional LMS where instantaneous gradient is used for weight update, in mLMS the momentum of the gradient change is used \cite{shynk1988lms}.  By incorporating the momentum of the gradient, the mLMS can achieve better convergence without compromising the steady-state error\cite{shynk1988lms, sharma1998analysis}.  Another method proposed by Raja et al \cite{LMS12}, is the application of fractional calculus.  Using the same method a series of papers have been published \cite{dubey2012flms, LMS17, LMS13}; claiming the improved steady-state and convergence performance.  But all these variants \cite{LMS12, LMS13, LMS16} have been criticized and are shown to offer no improvement over the conventional LMS \cite{bershad2017comments, ASLAM2015279}.  
Recently, another method in the same direction is proposed named as momentum FLMS (mFLMS) \cite{zubair2018momentum}. We argue that the proposed method does not \cite{zubair2018momentum} improve the performance of the conventional LMS and has serious flaws in its analysis presented in the paper.  Our argument is supported by analytical reasoning and extensive simulations.  

The organization of this paper is as follows: Flaws in the design, analysis and simulation setup of the mFLMS paper are discussed in section \ref{Sec:Issues}.  Estimation model, simulation setup and evaluation parameters are defined in  section \ref{Sec:Simulation}, followed by results and discussion in section \ref{Sec:RandD}.  Finally, the paper is concluded in section \ref{Sec:Con}.

\section{Remarks on \textquotedblleft Momentum fractional LMS for power signal parameter estimation \textquotedblright}\label{Sec:Issues}
\subsection{Design}
This section focuses on the main flaws of the design of the mFLMS algorithm \cite{zubair2018momentum}.  The structure of the mFLMS algorithm follows the architecture of the FLMS algorithm \cite{LMS12}, thus the problems of the FLMS that are mentioned in \cite{bershad2017comments}, are inherited in the mFLMS algorithm as well.  Let us first rewrite the expressions presented in \cite{zubair2018momentum}, they will be referred in the forthcoming mathematical analysis.  The first equation is the weight update rule for the FLMS algorithm presented in \cite{zubair2018momentum}:
\begin{equation}
\mathbf{\hat{w}}(n + 1) = \mathbf{\hat{w}}(n) + \mu_1 e(n) \textbf{u}(n) + \frac{\mu_f}{\Gamma(2 - f)} e(n) \textbf{u}(n) \odot \lvert\mathbf{\hat{w}}\rvert^{1 - f}(n) . \tag{12}\label{eq:12}
\end{equation}
where $ \mu_1 $ and $ \mu_f $ are the real positive values, defining the step-size.  While $ f $ is the fractional power in the range $ 0 < f < 1 $ and $\Gamma$ represents the Gamma function.  Here, $\textbf{u}$ and $\mathbf{\hat{w}}$ are the vectors defining the input and the estimated weights of the filter respectively.  The error between target and estimated outputs is $e(n)$.

The key set of equations derived in \cite{zubair2018momentum} are given as:
\begin{align}
& \mathbf{\hat{w}}(n + 1) = \mathbf{\hat{w}}(n) + \textbf{v}(n + 1), \tag{13}\label{eq:13}
\\
&\textbf{v}(n + 1) = \alpha \textbf{v}(n) + \textbf{g}(n), \tag{14}\label{eq:14}
\\
&\textbf{g}(n) = \mu_1 e(n) \textbf{u}(n) + \frac{\mu_f}{\Gamma(2 - f)} e(n) \textbf{u}(n) \odot \lvert\mathbf{\hat{w}}\rvert^{1 - f}(n), \tag{15}\label{eq:15}
\end{align}
where $\textbf{v}(n + 1)$ is the velocity update term defined by the momentum, $\alpha \in (0,1)$, and the instantaneous gradients $ \textbf{g}(n) $.  

According to \cite{zubair2018momentum}, Eq. \eqref{eq:12} is the weight update equation of the Fractional LMS algorithms \cite{LMS12}, which uses \textbf{\textit{absolute}} value of the weight vector primarily to avoid complex values.  However with this change, the fractional gradient term will not be the actual gradient of the cost function.  
Indeed, to exploit the fractional gradient characteristics, the complex domain information needs to be processed accordingly and proper gradient information must be incorporated \cite{CLMS}.  With this change of \textit{absolute} value of the weight vector, the FLMS algorithm has shown to offer no improvement compared to the conventional method \cite{bershad2017comments}.  Since the mFLMS approach \cite{zubair2018momentum} also utilizes the absolute values of $\mathbf{\hat{w}}$ in equation \eqref{eq:15} avoiding the complex mathematics, the same argument is valid for this case as well.

\subsection{Convergence Analysis}
%
There are serious flaws in the analysis section of the mFLMS algorithm, which are listed below.
\begin{enumerate} 
	\item  By employing the assumption of $ \mu_f = \mu_1 \Gamma(2 - f) $ , the authors construct the equation \eqref{eq:16} using equation \eqref{eq:13} to \eqref{eq:15}, as:
	\begin{equation}
	\mathbf{\hat{w}}(n + 1) = \mathbf{\hat{w}}(n) + \alpha \left[ \mathbf{\hat{w}}(n) - \mathbf{\hat{w}}(n - 1)\right]  + \mu_1 e(n) \textbf{u}(n) \odot \lvert\mathbf{\hat{w}}\rvert^{1 - f}(n) . \tag{16}\label{eq:16}
	\end{equation}
	However, the solved equation should be  
	\begin{equation*}
	\mathbf{\hat{w}}(n + 1) = \mathbf{\hat{w}}(n) + \alpha \left[ \mathbf{\hat{w}}(n) - \mathbf{\hat{w}}(n - 1)\right]  + \mu_1 e(n) \left[ \textbf{u}(n) + \textbf{u}(n) \odot \lvert\mathbf{\hat{w}}(n)\rvert^{1 - f}\right] . 
	\end{equation*}
	
	\item  Even if we consider the above mentioned flaw as a typo error, the solution from \textbf{(16)} to \textbf{(33)}, are not consistent and cannot be proved under any condition. 
	\item The last term in Eq. (17) should have $\odot$ operator (element by element operation), as it was used in Eqs. (15) and (16) of the paper.
	\item The last term in Eq. (17) should not have one added to the last bracket, rather the equation should be
	\begin{multline*}
	\mathbf{\Delta\hat{w}}(n + 1) = \mathbf{\Delta\hat{w}}(n) + \alpha \left[ \mathbf{\Delta\hat{w}}(n) - \mathbf{\Delta\hat{w}}(n - 1)\right]  \\+ \mu_1 e(n) \left[ \textbf{u}(n) + \textbf{u}(n) \odot \lvert\mathbf{w_{opt}} + \mathbf{\Delta\hat{w}}(n)\rvert^{1 - f}\right] .
	\end{multline*}
	

	\item Eq. (19) is incorrect because it defines the binomial formula for vectors.  The right hand side of the equation is however resulting in a scalar term (note that all the vectors are considered to be column vectors in this analysis). 
	\item The result in Eq. (23) is technically incorrect, in general the expectation of fractional power of any random variable cannot be replaced by the expectation of its unit power and a linear operation \cite{IDCF}. Thus, the expression in Eq. (23) is an approximation. 
	\item The convergence analysis of the mFLMS algorithm is primarily dependent on the evaluation of the expectation terms of the form $E[\Delta w^\gamma (n)]$ (where $\gamma$ is a fractional number). However, the analysis provided in this paper does not provide its accurate evaluation.  Instead, the authors have approximated this fractional moment using some function $G$ without providing its expression (see Eq. (23)). Thus, the whole convergence analysis is vague and the stability bound derived in Eq. (33) is meaningless without the knowledge of the function $G$.  
\end{enumerate}

\subsection{Simulation Setup}
In adaptive signal processing, performance comparison between algorithms can be made on the basis of different criteria.  Three important measures of performance are: (1) Convergence rate, (2) steady-state error and (3) computational complexity.  In \cite[Sec. 3.1]{zubair2018momentum}, it is already mentioned that the mFLMS algorithm is computationally very expensive, we therefore focus only on convergence and steady-state measures for our experiments.

We argue that in \cite{zubair2018momentum} the simulation parameters used for LMS algorithm are not appropriate.  With the learning rate of $1\times10^{-3}$ (as adopted in \cite{zubair2018momentum}) the LMS algorihtm converges slowly, interestingly the graphs are shown for just $4000$ iterations (see \cite[Fig. 2]{zubair2018momentum}).  For a fair evaluation, both algorithms must be setup at either equal convergence (for steady-state performance comparison) or equal steady-state (for convergence performance).  Also if one algorithm can perform better than the other in both aspects, then higher convergence rate at the cost of low steady-state performance must be shown. 
 
\section{Experimental Setup}\label{Sec:Simulation}
To re-evaluate the performance of the mFLMS and the conventional LMS algorithm, we considered the same problem of power signal parameter estimation as used in \cite{zubair2018momentum}.
\subsection{Power signal parameter estimation model}
In this section an overview of the power signal estimation model is provided. To estimate signal parameter, a sampled multi-harmonic sinusoidal signal $y$, with different amplitudes and phases, is considered, i.e.,
\begin{equation}\label{eq:multiharmonic}
	y (n) = \sum\limits_{k = 1}^{N} a_k \sin(n\omega_k + \phi_k ) + \epsilon(n) ,
\end{equation}
where $\epsilon(n) $ is a Gaussian disturbance of zero mean and constant variance $ \sigma^2 $, $a_k$, $\omega_k$ and $\phi_k$ are the amplitude, the angular frequency and the phase shifts of the sinusoid respectively.

With the help of the trigonometric identities, Eq. \eqref{eq:multiharmonic} can be transformed into:
\begin{equation}\label{eq:trigonometric}
	y (n) = \sum\limits_{k = 1}^{N} a_k \left(\sin n\omega_k \cos\phi_k + \cos n\omega_k \sin\phi_k\right)  + \epsilon(n) ,
\end{equation}
following the assumption made in \cite{zubair2018momentum}, the frequencies $\omega_k$ of the four sinusoids are taken to be known, Now, Eq. \eqref{eq:trigonometric} can be written as:
\begin{equation}\label{eq:frequency}
	y (n) = \sum\limits_{k = 1}^{N} b_k \sin n\omega_k  + c_k\cos n\omega_k   + \epsilon(n) ,
\end{equation}
where $b_k = a_k \cos\phi_k $ and $ c_k = a_k \sin\phi_k $ are the unknown parameters, $ a_k $ and $ \phi_k $ can be obtained using the following relations:
\begin{equation}
	a_k = \sqrt{b_k^2 + c_k^2}	\hspace{1cm}	\phi_k = \tan^{-1}\frac{c_k}{b_k} .
\end{equation}
The desired vector of parameters $ \mathbf{\theta}$, and the corresponding input vector are defined as:
\begin{equation}\label{eq:theta}
	\mathbf{\theta} = [b_1 , c_1 , b_2 , c_2 , \cdots , b_N , c_N]^T \in \rm I\!R^{2N} ,
\end{equation}
\begin{equation}\label{eq:phi}
	\mathbold{\psi} = [\sin \omega_1 n, \cos \omega_1 n, \sin \omega_2 n, \cos \omega_2 n, \cdots , \sin \omega_N n, \cos \omega_N n]^T \in \rm I\!R^{2N} .
\end{equation}

Finally the power signal parameter estimation model is given as:
\begin{equation}\label{eq:ID_Model}
	y (n) = \mathbold{\psi}^T(n) \theta + \epsilon(n) ,
\end{equation}

While applying the LMS and mFLMS algorithms, the unknown parameter $ \theta $ will be treated as the weight vector and will be updated using the respective weight update rules.

\subsection{Simulation setup}\label{Sec:4}
We consider a composite signal of sinusoidal with four different frequencies \cite{xu2017recursive}, i.e.,
\begin{multline}\label{signal_f}
	y (n) = 1.8 \sin (0.07 n + 0.95) + 2.9 \sin (0.5 n + 0.8) \\ + 4 \sin (2 n + 0.76) + 2.5 \sin (1.6 n + 1.1) + \epsilon(n) .
\end{multline}
We assume that all four frequencies used in equation (\ref{signal_f}) are known.  This assumption is the same made in \cite{zubair2018momentum}.  The desired parameter set is defined as:
\begin{equation}
	\theta = [ a_1 , a_2 , a_3 , a_4 , \phi_1 , \phi_2 , \phi_3 , \phi_4]^T = [1.8, 2.9, 4, 2.5, 0.95, 0.8, 0.76, 1.1]^T .
\end{equation}

\subsection{Evaluation metrics}\label{Sec:4.1}
The performance of both the LMS and the mFLMS algorithm is evaluted on two metrics: (1) mean squared error (MSE) and (2) the normalized weight differences (NWD) defined as a fitness function $\delta$ in \cite{zubair2018momentum}.  The NWD is given as:

\begin{equation}
	{\rm \mathbf{NWD}}(n) = \frac{\lvert\lvert \mathbold{\hat{\theta}(n)} - \mathbold{\theta} \rvert\rvert}{\lvert\lvert \mathbold{\theta} \rvert\rvert} ,
\end{equation}
where $\mathbold{\theta} $ and $ \hat{\mathbold{\theta}}$ are the desired and approximated parameter vectors at the $n$th iteration.  Another performance measure based on MSE is given as:
\begin{equation}
	{\rm MSE} = \frac{1}{M} \sum\limits_{i = 1}^{M} (\mathbold{\theta}_i - \mathbold{\hat{\theta}}_i)^2 ,
\end{equation}
where $ M $ is the length of the desired vector $\mathbold{\theta}$, and $ \hat{\mathbold{\theta}}$ is the final estimated vector of parameters.

\section{Results and Discussion}\label{Sec:RandD}
\begin{figure}
\centering
\begin{tabular}[b]{c}
	\centerline{\includegraphics[width=0.60\linewidth]{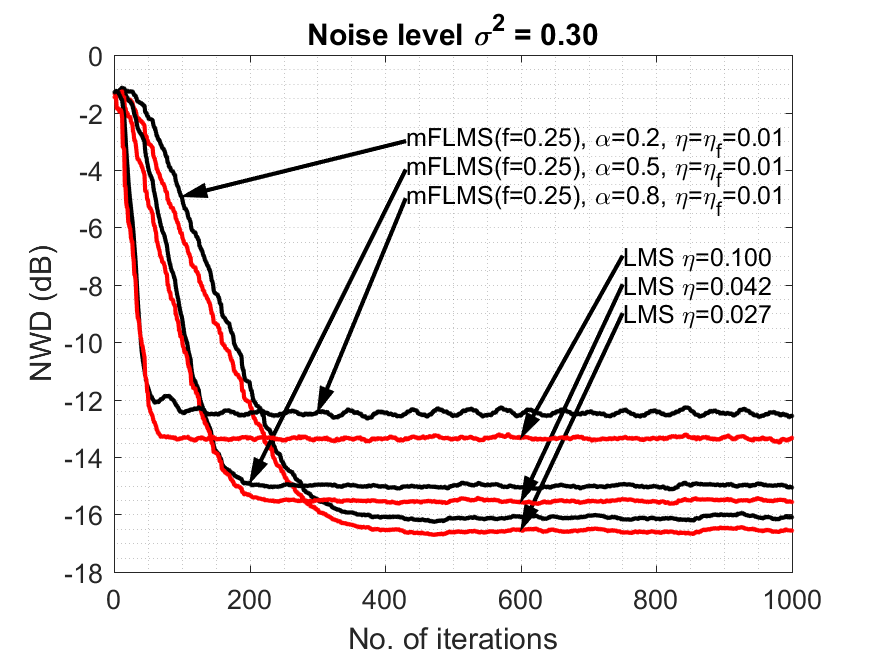}}
\\
	\small (a) $f=0.25$.
\end{tabular} \qquad
\begin{tabular}[b]{c}
	\centerline{\includegraphics[width=0.60\linewidth]{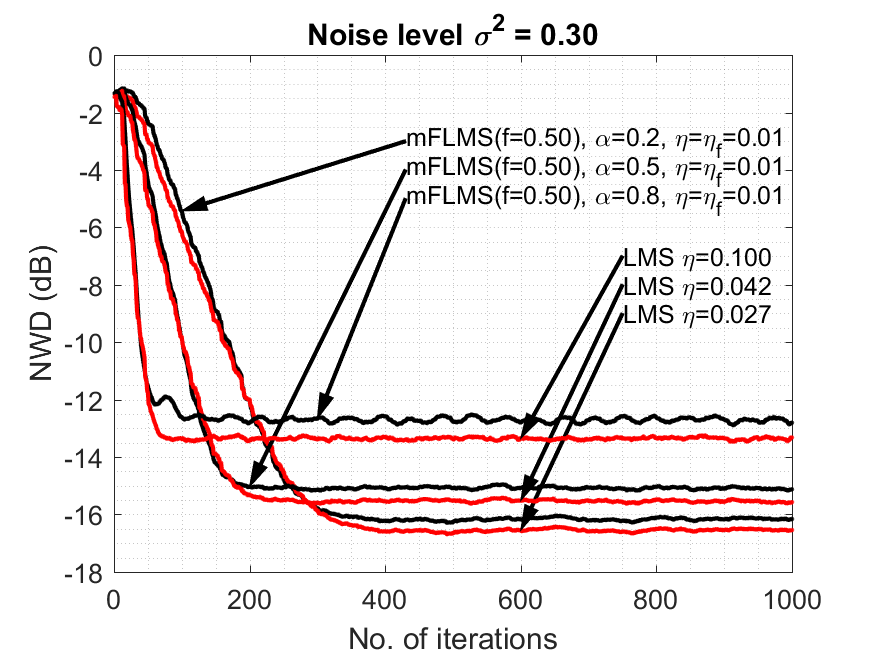}}
\\
	\small (b) $f=0.50$.
\end{tabular}
\begin{tabular}[b]{c}
		\centerline{\includegraphics[width=0.60\linewidth]{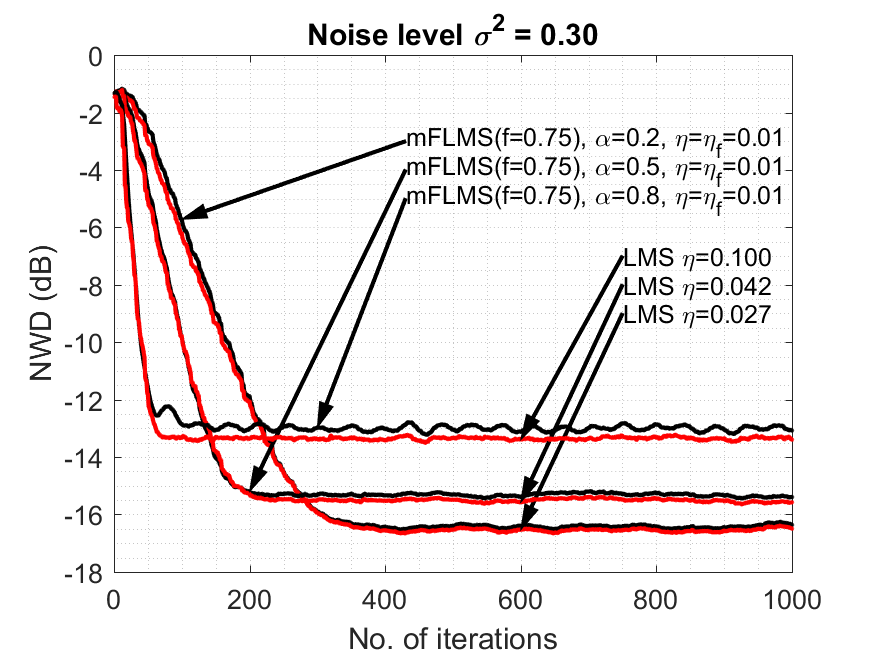}}
		 \\
	\small (c) $f=0.75$.
\end{tabular}
\caption{Comparison of steady state results for $\sigma^2 = 0.30$.}
\label{nl1}
\end{figure}

\begin{figure}
	\centering
	\begin{tabular}[b]{c}
			\centerline{\includegraphics[width=0.60\linewidth]{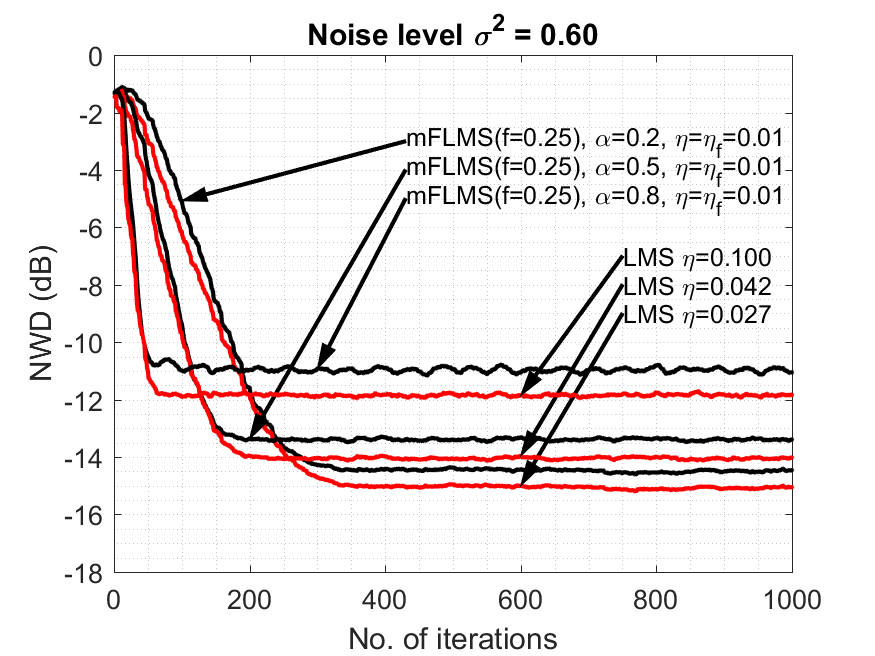}}
		\\
		\small (a) $f=0.25$.
	\end{tabular} \qquad
	\begin{tabular}[b]{c}
			\centerline{\includegraphics[width=0.60\linewidth]{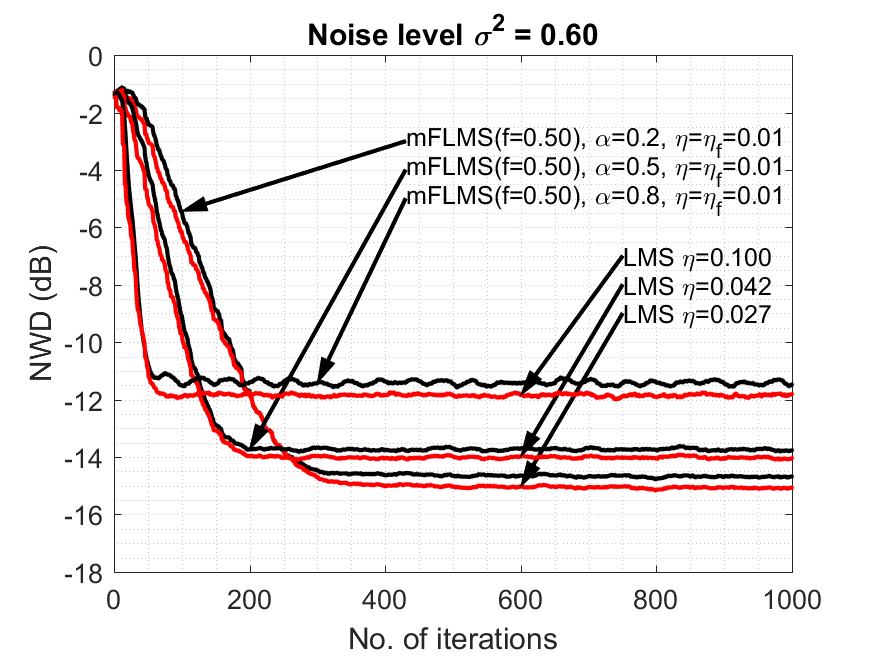}}
		 \\
		\small (b) $f=0.50$.
	\end{tabular}
	\begin{tabular}[b]{c}
			\centerline{\includegraphics[width=0.60\linewidth]{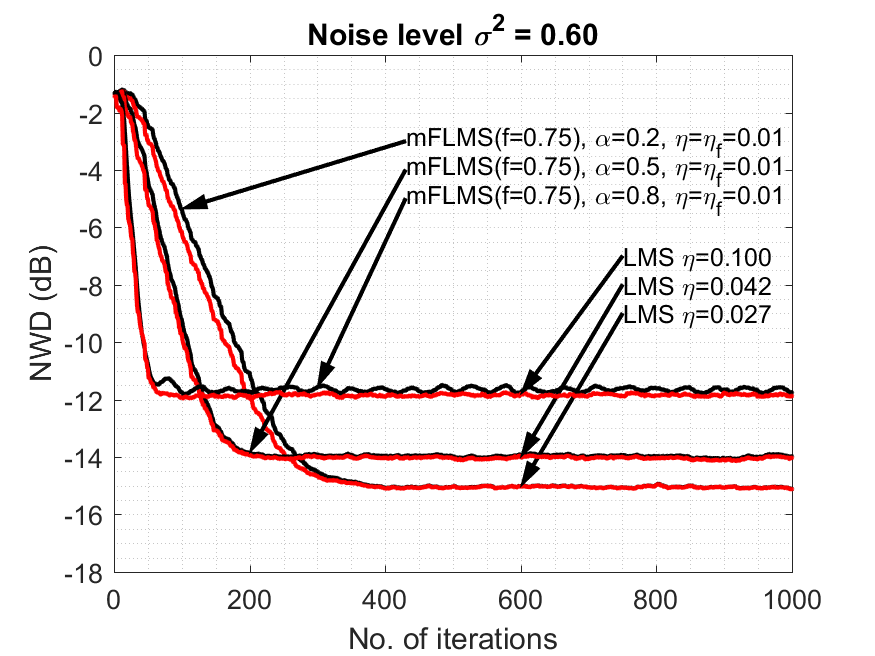}}
		 \\
		\small (c) $f=0.75$.
	\end{tabular}
\caption{Comparison of steady state results for $\sigma^2 = 0.60$.}
	\label{nl2}
\end{figure}

\begin{figure}
	\centering
	\begin{tabular}[b]{c}
			\centerline{\includegraphics[width=0.60\linewidth]{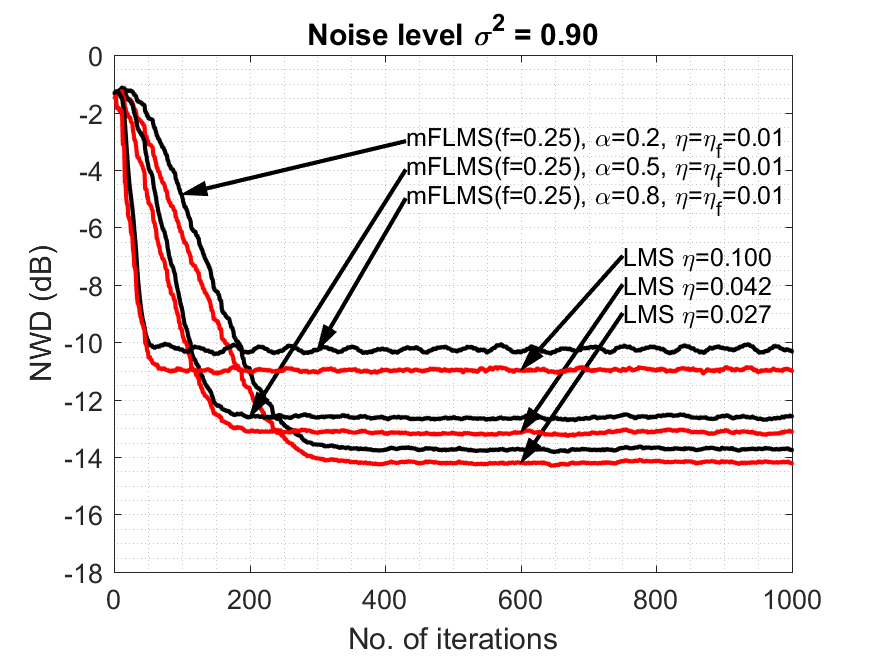}}
		\\
		\small (a) $f=0.25$.
	\end{tabular} \qquad
	\begin{tabular}[b]{c}
			\centerline{\includegraphics[width=0.60\linewidth]{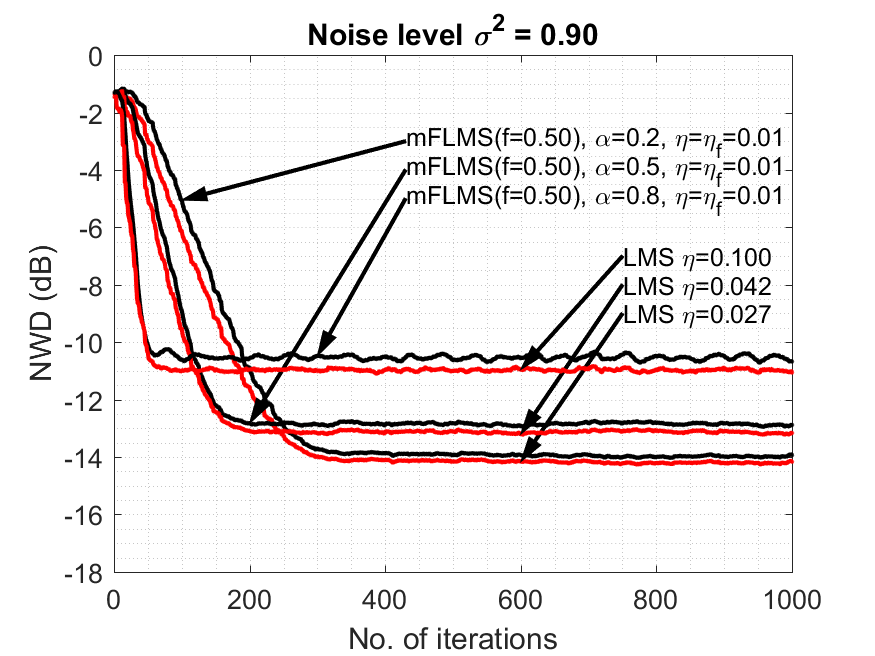}}
			 \\
		\small (b) $f=0.50$.
	\end{tabular}
	\begin{tabular}[b]{c}
			\centerline{\includegraphics[width=0.60\linewidth]{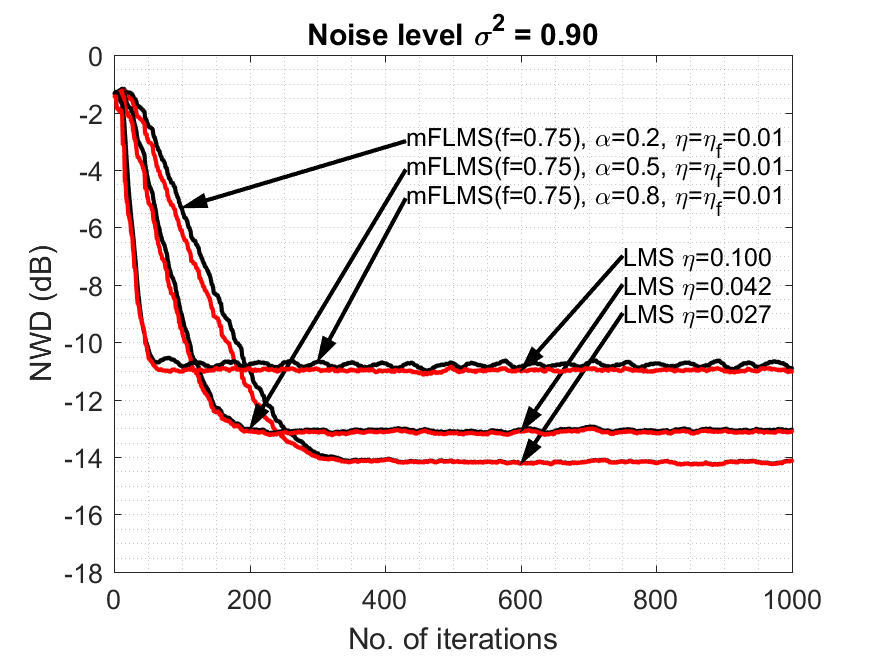}}
			 \\
		\small (c) $f=0.75$.
	\end{tabular}
	\caption{Comparison of steady state results for $\sigma^2 = 0.90$.}
	\label{nl3}
\end{figure}

The LMS and mFLMS \cite{zubair2018momentum} algorithm are compared for three noise levels, i.e., ($\sigma^2=0.3$, $\sigma^2=0.6$, and $\sigma^2=0.9$) .  Figs. \ref{nl1}, \ref{nl2}, and \ref{nl3} show the learning curves for mFLMS and LMS for respective noise levels.  We setup the LMS and the mFLMS algorithms for the same convergence performance, and compared the steady-state performance of the two.  From the Figs. \ref{nl1},\ref{nl2} and \ref{nl3}, it can be seen that under all conditions the LMS algorithm performed better than the mFLMS.  The experiments were repeated for $1000$ independent rounds and mean results for every $100$th iteration was reported.  For each independent round, the weights were initialized using random values obtained from Gaussian distribution of mean zero and variance of one.  The performance of the mFLMS algorithm is evaluated for three different values of momentum $\alpha$, i.e., ($\alpha = 0.2, 0.5, 0.8$) and for each value of the momentum, the algorithm is compared for three different fractional powers, i.e., $0.25$, $0.5$, and $0.75$.  We observed that the effect of momentum in the mFLMS algorithm is equivalent to that of the learning rate in the conventional LMS, i.e, the increase in convergence rate increases the steady-state error.  The effect of fractional power $f$ is also arguable, it is shown to degrade the steady-state performance without even improving the convergence rate.  Normalized weight difference results for every $100$th iteration under all $27$ scenarios ($3$ noise level, $3$ different step-sizes, and $3$ different values of fractional power) are summarized in Tables \ref{nlt1}, \ref{nlt2} and \ref{nlt3}.  Final obtained values of estimation parameters and mean squared error (MSE) for noise levels ($\sigma^2=0.3$, $\sigma^2=0.6$, and $\sigma^2=0.9$) are reported in Table \ref{nlm1}, \ref{nlm2} and \ref{nlm3} respectively.

\newpage

\begin{table}[h!]
	\tiny
	\centering
	\caption{Performance comparison based on fitness achieved at specific iterations for $\sigma^{2} = 0.30^{2}$}
	\label{nlt1}
	\begin{tabular}{@{}llllllllllll@{}}
		\toprule
		\textbf{Method}            & \textbf{$\alpha$} & \multicolumn{10}{l}{\textbf{Fitness achieved at specific iterations}}                   \\ \midrule
		\textbf{}                  & \textbf{}         & 100    & 200    & 300    & 400    & 500    & 600    & 700    & 800    & 900    & 1000   \\
		\textbf{mFLMS($f$=0.25)}     & \textbf{0.2}      & 0.3218 & 0.0725 & 0.0284 & 0.0246 & 0.0243 & 0.0246 & 0.0250 & 0.0242 & 0.0250 & 0.0246 \\
		\textbf{mFLMS($f$=0.50)}     & \textbf{0.2}      & 0.2867 & 0.0628 & 0.0264 & 0.0239 & 0.0237 & 0.0242 & 0.0239 & 0.0245 & 0.0242 & 0.0242 \\
		\textbf{mFLMS($f$=0.75)}     & \textbf{0.2}      & 0.2679 & 0.0647 & 0.0257 & 0.0228 & 0.0229 & 0.0229 & 0.0229 & 0.0228 & 0.0227 & 0.0231 \\
		\textbf{LMS($\eta$=0.027)} &                   & \textbf{0.2397} & \textbf{0.0619} & \textbf{0.0260} & \textbf{0.0222} & \textbf{0.0218} & \textbf{0.0222} & \textbf{0.0224} & \textbf{0.0217} & \textbf{0.0225} & \textbf{0.0221} \\
		\textbf{mFLMS($f$=0.25)}     & \textbf{0.5}      & 0.1228 & 0.0323 & 0.0319 & 0.0316 & 0.0316 & 0.0312 & 0.0316 & 0.0314 & 0.0317 & 0.0312 \\
		\textbf{mFLMS($f$=0.50)}     & \textbf{0.5}      & 0.1033 & 0.0314 & 0.0306 & 0.0311 & 0.0310 & 0.0311 & 0.0312 & 0.0313 & 0.0308 & 0.0307 \\
		\textbf{mFLMS($f$=0.75)}     & \textbf{0.5}      & 0.1044 & 0.0302 & 0.0294 & 0.0293 & 0.0295 & 0.0289 & 0.0302 & 0.0298 & 0.0292 & 0.0289 \\
		\textbf{LMS($\eta$=0.042)} &                   & \textbf{0.1041} & \textbf{0.0294} & \textbf{0.0281} & \textbf{0.0279} & \textbf{0.0281} & \textbf{0.0278} & \textbf{0.0281} & \textbf{0.0279} & \textbf{0.0281} & \textbf{0.0278} \\
		\textbf{mFLMS($f$=0.25)}     & \textbf{0.8}      & 0.0572 & 0.0567 & 0.0571 & 0.0569 & 0.0553 & 0.0556 & 0.0571 & 0.0572 & 0.0564 & 0.0554 \\
		\textbf{mFLMS($f$=0.50)}     & \textbf{0.8}      & 0.0548 & 0.0525 & 0.0543 & 0.0553 & 0.0529 & 0.0525 & 0.0535 & 0.0546 & 0.0537 & 0.0524 \\
		\textbf{mFLMS($f$=0.75)}     & \textbf{0.8}      & 0.0518 & 0.0500 & 0.0505 & 0.0499 & 0.0495 & 0.0488 & 0.0502 & 0.0516 & 0.0497 & 0.0492 \\
		\textbf{LMS($\eta$=0.1)}    &                   & \textbf{0.0462} & \textbf{0.0466} & \textbf{0.0462} & \textbf{0.0462} & \textbf{0.0461} & \textbf{0.0464} & \textbf{0.0459} & \textbf{0.0460} & \textbf{0.0466} & \textbf{0.0464} \\
		\bottomrule
	\end{tabular}
\end{table}

\begin{table}[h!]
	\tiny
	\centering
	\caption{Performance comparison based on fitness achieved at specific iterations for $\sigma^{2} = 0.60^{2}$}
	\label{nlt2}
	\begin{tabular}{@{}llllllllllll@{}}
		\toprule
		\textbf{Method}            & \textbf{$\alpha$} & \multicolumn{10}{l}{\textbf{Fitness achieved at specific iterations}}                   \\ \midrule
		\textbf{}                  & \textbf{}         & 100    & 200    & 300    & 400    & 500    & 600    & 700    & 800    & 900    & 1000   \\
		\textbf{mFLMS($f$=0.25)}     & \textbf{0.2}      & 0.3108 & 0.0694 & 0.0373 & 0.0362 & 0.0366 & 0.0361 & 0.0359 & 0.0351 & 0.0355 & 0.0358 \\
		\textbf{mFLMS($f$=0.50)}     & \textbf{0.2}      & 0.2857 & 0.0676 & 0.0357 & 0.0344 & 0.0342 & 0.0344 & 0.0341 & 0.0334 & 0.0339 & 0.0340 \\
		\textbf{mFLMS($f$=0.75)}     & \textbf{0.2}      & 0.2912 & 0.0782 & 0.0354 & 0.0311 & 0.0314 & 0.0315 & 0.0313 & 0.0318 & 0.0311& 0.0308 \\
		\textbf{LMS($\eta$=0.027)} &                   & \textbf{0.2428} & \textbf{0.0656} & \textbf{0.0338} & \textbf{0.0315} & \textbf{0.0320} & \textbf{0.0316} & \textbf{0.0313} & \textbf{0.0307} & \textbf{0.0309} & \textbf{0.0307} \\
		\textbf{mFLMS($f$=0.25)}     & \textbf{0.5}      & 0.1158 & 0.0460 & 0.0459 & 0.0467 & 0.0457 & 0.0466 & 0.0458 & 0.0462 & 0.0459 & 0.0458 \\
		\textbf{mFLMS($f$=0.50)}     & \textbf{0.5}      & 0.1235 & 0.0424 & 0.0430 & 0.0421 & 0.0426 & 0.0426 & 0.0427 & 0.0420 & 0.0419 & 0.0420 \\
		\textbf{mFLMS($f$=0.75)}     & \textbf{0.5}      & 0.1198 & 0.0412 & 0.0406 & 0.0402 & 0.0400 & 0.0405 & 0.0398 & 0.0400 & 0.0394 & 0.0400 \\
		\textbf{LMS($\eta$=0.042)} &                   & \textbf{0.1077} & \textbf{0.0402} & \textbf{0.0394} & \textbf{0.0400} & \textbf{0.0392} & \textbf{0.0402} & \textbf{0.0395} & \textbf{0.0398} & \textbf{0.0392} & \textbf{0.0396} \\
		\textbf{mFLMS($f$=0.25)}     & \textbf{0.8}      & 0.0797 & 0.0801 & 0.0807 & 0.0815 & 0.0804 & 0.0784 & 0.0801 & 0.0829 & 0.0801 & 0.0785 \\
		\textbf{mFLMS($f$=0.50)}     & \textbf{0.8}      & 0.0708 & 0.0724 & 0.0722 & 0.0731 & 0.0716 & 0.0718 & 0.0739 & 0.0724 & 0.0719 & 0.0712 \\
		\textbf{mFLMS($f$=0.75)}     & \textbf{0.8}      & 0.0669 & 0.0672 & 0.0697 & 0.0687 & 0.0676 & 0.0669 & 0.0680 & 0.0690 & 0.0684 & 0.0668 \\
		\textbf{LMS($\eta$=0.1)}    &                   & \textbf{0.0646} & \textbf{0.0653} & \textbf{0.0647} & \textbf{0.0657} & \textbf{0.0664} & \textbf{0.0652} & \textbf{0.0649} & \textbf{0.0665} & \textbf{0.0656} & \textbf{0.0654} \\
		\bottomrule
	\end{tabular}
\end{table}

\begin{table}[h!]
	\tiny
	\centering
	\caption{Performance comparison based on fitness achieved at specific iterations for $\sigma^{2} = 0.90^{2}$}
	\label{nlt3}
	\begin{tabular}{@{}llllllllllll@{}}
		\toprule
		\textbf{Method}            & \textbf{$\alpha$} & \multicolumn{10}{l}{\textbf{Fitness achieved at specific iterations}}                   \\ \midrule
		\textbf{}                  & \textbf{}         & 100    & 200    & 300    & 400    & 500    & 600    & 700    & 800    & 900    & 1000   \\
		\textbf{mFLMS($f$=0.25)}     & \textbf{0.2}      & 0.3273 & 0.0811 & 0.0439 & 0.0422 & 0.0422 & 0.0422 & 0.0423 & 0.0430 & 0.0426 & 0.0420 \\
		\textbf{mFLMS($f$=0.50)}     & \textbf{0.2}      & 0.3125 & 0.0815 & 0.0426 & 0.0411 & 0.0405 & 0.0404 & 0.0405 & 0.0398 & 0.0400 & 0.0404 \\
		\textbf{mFLMS($f$=0.75)}     & \textbf{0.2}      & 0.2941 & 0.0816 & 0.0410 & 0.0390 & 0.0384 & 0.0378 & 0.0385 & 0.0380 & 0.0378 & 0.0387 \\
		\textbf{LMS($\eta$=0.027)} &                   & \textbf{0.2423} & \textbf{0.0693} & \textbf{0.0393} & \textbf{0.0378} & \textbf{0.0378} & \textbf{0.0380} & \textbf{0.0380} & \textbf{0.0386} & \textbf{0.0382} & \textbf{0.0378} \\
		\textbf{mFLMS($f$=0.25)}     & \textbf{0.5}      & 0.1319 & 0.0549 & 0.0547 & 0.0544 & 0.0546 & 0.0543 & 0.0543 & 0.0551 & 0.0551 & 0.0552 \\
		\textbf{mFLMS($f$=0.50)}     & \textbf{0.5}      & 0.1288 & 0.0522 & 0.0518 & 0.0516 & 0.0522 & 0.0510 & 0.0523 & 0.0523 & 0.0521 & 0.0515 \\
		\textbf{mFLMS($f$=0.75)}     & \textbf{0.5}      & 0.1248 & 0.0498 & 0.0496 & 0.0495 & 0.0489 & 0.0496 & 0.0507 & 0.0489 & 0.0497 & 0.0495 \\
		\textbf{LMS($\eta$=0.042)} &                   & \textbf{0.1110} & \textbf{0.0488} & \textbf{0.0483} & \textbf{0.0481} & \textbf{0.0483} & \textbf{0.0485} & \textbf{0.0480} & \textbf{0.0490} & \textbf{0.0489} & \textbf{0.0488} \\
		\textbf{mFLMS($f$=0.25)}     & \textbf{0.8}      & 0.0944 & 0.0925 & 0.0939 & 0.0945 & 0.0933 & 0.0922 & 0.0956 & 0.0955 & 0.0951 & 0.0928 \\
		\textbf{mFLMS($f$=0.50)}     & \textbf{0.8}      & 0.0869 & 0.0880 & 0.0896 & 0.0892 & 0.0867 & 0.0880 & 0.0899 & 0.0893 & 0.0894 & 0.0856 \\
		\textbf{mFLMS($f$=0.75)}     & \textbf{0.8}      & 0.0826 & 0.0837 & 0.0854 & 0.0844 & 0.0844 & 0.0821 & 0.0849 & 0.0843 & 0.0832 & 0.0811 \\
		\textbf{LMS($\eta$=0.1)}    &                   & \textbf{0.0800} & \textbf{0.0794} & \textbf{0.0786} & \textbf{0.0801} & \textbf{0.0798} & \textbf{0.0799} & \textbf{0.0806} & \textbf{0.0797} & \textbf{0.0809} & \textbf{0.0794} \\
		\bottomrule
	\end{tabular}
\end{table}
\newpage

\begin{table}[h]
	\tiny
	\centering
	\caption{Performance comparison based on estimation error for $\sigma^{2} = 0.30^{2}$}
	\label{nlm1}
	\begin{tabular}{lcccccccccc}
		\hline
		\textbf{Method}            & \textbf{$\alpha$}             & \multicolumn{8}{l}{\textbf{Adaptive parameters}}                                                                                              & \textbf{MSE}    \\ \hline
		\textbf{}                  & \multicolumn{1}{l}{\textbf{}} & $\theta_{1}$    & $\theta_{2}$    & $\theta_{3}$    & $\theta_{4}$    & $\theta_{5}$    & $\theta_{6}$    & $\theta_{7}$    & $\theta_{8}$    &                 \\
		\textbf{mFLMS($f$=0.25)} & \textbf{0.2}                  & 1.9105          & 2.8816          & 3.8812          & 2.5303          & 0.9977          & 0.8001          & 0.7639          & 1.1139          & 3.98E-06           \\
		\textbf{mFLMS($f$=0.50)} & \textbf{0.2}                  & 1.8005          & 2.9045          & 3.9970          & 2.5002          & 0.9506          & 0.7991          & 0.7598          & 1.1010          & 3.95E-06          \\
		\textbf{mFLMS($f$=0.75)} & \textbf{0.2}                  & 1.8006          & 2.9043          & 4.0026          & 2.5023          & 0.9496          & 0.8003          & 0.7599          & 1.1008          & 7.71E-07          \\
		\textbf{LMS($\eta$=0.027)} & \multicolumn{1}{l}{\textbf{}} & \textbf{1.8017} & \textbf{2.9048} & \textbf{3.9983} & \textbf{2.5009} & \textbf{0.9505} & \textbf{0.7992} & \textbf{0.7598} & \textbf{1.1008} & \textbf{6.41E-07} \\
		\textbf{mFLMS($f$=0.25)} & \textbf{0.5}                  & 1.7999          & 2.8985          & 3.9991          & 2.4947          & 0.9500          & 0.7985          & 0.7583          & 1.1004          & 9.19E-06          \\
		\textbf{mFLMS($f$=0.50)} & \textbf{0.5}                  & 1.8005          & 2.9028          & 3.9954          & 2.5000          & 0.9519          & 0.8016          & 0.7602          & 1.1000          & 4.52E-06          \\
		\textbf{mFLMS($f$=0.75)} & \textbf{0.5}                  & 1.7986          & 2.9042          & 4.0059          & 2.5043          & 0.9501          & 0.8010          & 0.7605          & 1.1001          & 4.37E-06          \\
		\textbf{LMS($\eta$=0.042)} & \multicolumn{1}{l}{\textbf{}} & \textbf{1.8022} & \textbf{2.9034} & \textbf{3.9975} & \textbf{2.5013} & \textbf{0.9515} & \textbf{0.8013} & \textbf{0.7600} & \textbf{1.0999} & \textbf{1.90E-06} \\
		\textbf{mFLMS($f$=0.25)} & \textbf{0.8}                  & 1.8080          & 2.8936          & 3.9853          & 2.5001          & 0.9457          & 0.8011          & 0.7614          & 1.1020          & 4.33E-05          \\
		\textbf{mFLMS($f$=0.50)} & \textbf{0.8}                  & 1.8022          & 2.9084          & 4.0051          & 2.5023          & 0.9431          & 0.7997          & 0.7628          & 1.0994          & 3.92E-05          \\
		\textbf{mFLMS($f$=0.75)} & \textbf{0.8}                  & 1.8102          & 2.9010          & 4.0079          & 2.5113          & 0.9538          & 0.7985          & 0.7611          & 1.1006          & 2.03E-05          \\
		\textbf{LMS($\eta$=0.1)}   & \multicolumn{1}{l}{\textbf{}} & \textbf{1.8083} & \textbf{2.9036} & \textbf{4.0086} & \textbf{2.5127} & \textbf{0.9514} & \textbf{0.7982} & \textbf{0.7616} & \textbf{1.0994} & \textbf{1.51E-05} \\
		\textbf{True values}       & \multicolumn{1}{l}{\textbf{}} & \textbf{1.8}    & \textbf{2.9}    & \textbf{4}      & \textbf{2.5}    & \textbf{0.95}   & \textbf{0.8}    & \textbf{0.76}   & \textbf{1.1}    & \textbf{0}      \\ \hline
	\end{tabular}
\end{table}

\begin{table}[h!]
	\tiny
	\centering
	\caption{Performance comparison based on estimation error for $\sigma^{2} = 0.60^{2}$}
	\label{nlm2}
	\begin{tabular}{lcccccccccc}
		\hline
		\textbf{Method}            & \textbf{$\alpha$}             & \multicolumn{8}{l}{\textbf{Adaptive parameters}}                                                                                              & \textbf{MSE}    \\ \hline
		\textbf{}                  & \multicolumn{1}{l}{\textbf{}} & $\theta_{1}$    & $\theta_{2}$    & $\theta_{3}$    & $\theta_{4}$    & $\theta_{5}$    & $\theta_{6}$    & $\theta_{7}$    & $\theta_{8}$    &                 \\
		\textbf{mFLMS($f$=0.25)} & \textbf{0.2}                  & 1.7995          & 2.8981          & 4.0013          & 2.5030          & 0.9482          & 0.7999          & 0.7604          & 1.1011          & 3.94E-06          \\
		\textbf{mFLMS($f$=0.50)} & \textbf{0.2}                  & 1.7992          & 2.8993          & 4.0025          & 2.4998          & 0.9495          & 0.7993          & 0.7610          & 1.1022          & 2.43E-06          \\
		\textbf{mFLMS($f$=0.75)} & \textbf{0.2}                  & 1.8004          & 2.9036          & 4.0037          & 2.5012          & 0.9504          & 0.7985          & 0.7597          & 1.1003          & 9.96E-07          \\
		\textbf{LMS($\eta$=0.027)} & \multicolumn{1}{l}{\textbf{}} & \textbf{1.8925} & \textbf{2.9328} & \textbf{3.7577} & \textbf{2.4026} & \textbf{0.9750} & \textbf{0.8300} & \textbf{0.7450} & \textbf{1.0707} & \textbf{9.21E-07} \\
		\textbf{mFLMS($f$=0.25)} & \textbf{0.5}                  & 1.7981          & 2.8967          & 3.9930          & 2.4939          & 0.9537          & 0.7989          & 0.7579          & 1.1008          & 1.52E-05          \\
		\textbf{mFLMS($f$=0.50)} & \textbf{0.5}                  & 1.7993          & 2.9015          & 3.9998          & 2.5025          & 0.9512          & 0.8001          & 0.7595          & 1.0979          & 7.78E-06          \\
		\textbf{mFLMS($f$=0.75)} & \textbf{0.5}                  & 1.8063          & 2.9002          & 4.0020          & 2.4988          & 0.9477          & 0.7990          & 0.7621          & 1.1000          & 7.02E-06          \\
		\textbf{LMS($\eta$=0.042)} & \multicolumn{1}{l}{\textbf{}} & \textbf{1.8042} & \textbf{2.9021} & \textbf{3.9996} & \textbf{2.5007} & \textbf{0.9530} & \textbf{0.7994} & \textbf{0.7588} & \textbf{1.1003} & \textbf{4.20E-06} \\
		\textbf{mFLMS($f$=0.25)} & \textbf{0.8}                  & 1.8009          & 2.9186          & 3.9866          & 2.4943          & 0.9484          & 0.7994          & 0.7596          & 1.1051          & 7.35E-05         \\
		\textbf{mFLMS($f$=0.50)} & \textbf{0.8}                  & 1.8153          & 2.8892          & 3.9982          & 2.5124          & 0.9466          & 0.8009          & 0.7615          & 1.0984          & 6.53E-05          \\
		\textbf{mFLMS($f$=0.75)} & \textbf{0.8}                  & 1.8104          & 2.9010          & 4.0098          & 2.5102          & 0.9459          & 0.7992          & 0.7575          & 1.0927          & 4.81E-05          \\
		\textbf{LMS($\eta$=0.1)}   & \multicolumn{1}{l}{\textbf{}} & \textbf{1.8085} & \textbf{2.9238} & \textbf{3.9993} & \textbf{2.5115} & \textbf{0.9460} & \textbf{0.8009} & \textbf{0.7606} & \textbf{1.1027} & \textbf{3.58E-05} \\
		\textbf{True values}       & \multicolumn{1}{l}{\textbf{}} & \textbf{1.8}    & \textbf{2.9}    & \textbf{4}      & \textbf{2.5}    & \textbf{0.95}   & \textbf{0.8}    & \textbf{0.76}   & \textbf{1.1}    & \textbf{0}      \\ \hline
	\end{tabular}
\end{table}

\begin{table}[h!]
	\tiny
	\centering
	\caption{Performance comparison based on estimation error for $\sigma^{2} = 0.90^{2}$}
	\label{nlm3}
	\begin{tabular}{lcccccccccc}
		\hline
		\textbf{Method}            & \textbf{$\alpha$}             & \multicolumn{8}{l}{\textbf{Adaptive parameters}}                                                                                              & \textbf{MSE}    \\ \hline
		\textbf{}                  & \multicolumn{1}{l}{\textbf{}} & $\theta_{1}$    & $\theta_{2}$    & $\theta_{3}$    & $\theta_{4}$    & $\theta_{5}$    & $\theta_{6}$    & $\theta_{7}$    & $\theta_{8}$    &                 \\
		\textbf{mFLMS($f$=0.25)} & \textbf{0.2}                  & 1.7970          & 2.8958          & 3.9938          & 2.4978          & 0.9490          & 0.8003          & 0.7612          & 1.0989          & 1.59E-05          \\
		\textbf{mFLMS($f$=0.50)} & \textbf{0.2}                  & 1.8022          & 2.9038          & 3.9917          & 2.5003          & 0.9507          & 0.7986          & 0.7595          & 1.1001          & 1.14E-05          \\
		\textbf{mFLMS($f$=0.75)} & \textbf{0.2}                  & 1.8029          & 2.9100          & 3.9975          & 2.4966          & 0.9494          & 0.7998          & 0.7606          & 1.1010          & 9.29E-06          \\
		\textbf{LMS($\eta$=0.027)} & \multicolumn{1}{l}{\textbf{}} & \textbf{1.8025} & \textbf{2.9009} & \textbf{3.9986} & \textbf{2.5028} & \textbf{0.9486} & \textbf{0.8002} & \textbf{0.7612} & \textbf{1.0983} & \textbf{2.93E-06} \\
		\textbf{mFLMS($f$=0.25)} & \textbf{0.5}                  & 1.8015          & 2.8920          & 3.9919          & 2.4974          & 0.9508          & 0.7998          & 0.7603          & 1.1043          & 1.98E-05          \\
		\textbf{mFLMS($f$=0.50)} & \textbf{0.5}                  & 1.7986          & 2.9000          & 4.0045          & 2.4978          & 0.9489          & 0.7980          & 0.7600          & 1.0999          & 1.00E-05          \\
		\textbf{mFLMS($f$=0.75)} & \textbf{0.5}                  & 1.7968          & 2.8962          & 3.9955          & 2.5008          & 0.9520          & 0.7949          & 0.7620          & 1.0995          & 6.48E-06          \\
		\textbf{LMS($\eta$=0.042)} & \multicolumn{1}{l}{\textbf{}} & \textbf{1.8076} & \textbf{2.8993} & \textbf{3.9972} & \textbf{2.5052} & \textbf{0.9498} & \textbf{0.7997} & \textbf{0.7598} & \textbf{1.1035} & \textbf{4.02E-06} \\
		\textbf{mFLMS($f$=0.25)} & \textbf{0.8}                  & 1.8053          & 2.8980          & 3.9874          & 2.4987          & 0.9571          & 0.8026          & 0.7554          & 1.1052          & 1.33E-04          \\
		\textbf{mFLMS($f$=0.50)} & \textbf{0.8}                  & 1.8131          & 2.8926          & 3.9913          & 2.4938          & 0.9550          & 0.7984          & 0.7579          & 1.1071          & 5.28E-05          \\
		\textbf{mFLMS($f$=0.75)} & \textbf{0.8}                  & 1.8284          & 2.9032          & 3.9920          & 2.5063          & 0.9593          & 0.7984          & 0.7608          & 1.1070          & 3.72E-05          \\
		\textbf{LMS($\eta$=0.1)}   & \multicolumn{1}{l}{\textbf{}} & \textbf{1.8114} & \textbf{2.9164} & \textbf{4.0052} & \textbf{2.5240} & \textbf{0.9525} & \textbf{0.8016} & \textbf{0.7574} & \textbf{1.1014} & \textbf{3.32E-05} \\
		\textbf{True values}       & \multicolumn{1}{l}{\textbf{}} & \textbf{1.8}    & \textbf{2.9}    & \textbf{4}      & \textbf{2.5}    & \textbf{0.95}   & \textbf{0.8}    & \textbf{0.76}   & \textbf{1.1}    & \textbf{0}      \\ \hline
	\end{tabular}
\end{table}

\newpage
\section{Conclusion}\label{Sec:Con}
Recently, momemtum fractional LMS (mFLMS) algorithm has been proposed for the estimation of power signal parameter \cite{zubair2018momentum}.  The algorithm is shown to outperform the conventional LMS algorithm.  The mathematical assumptions made for a simplified derivation of the mFLMS algorithm are however invalid.  In this research, we have highlighted the discrepancies in the mFLMS algorithm proposed in \cite{zubair2018momentum}.  Extensive experiments have also been performed to thoroughly investigate the merits of the mFLMS algorithm.  After a careful analysis of the experimental results, we conclude that, under no condition the mFLMS algorithm is better than the conventional LMS algorithm for power signal parameter estimation.  In fact the LMS algorithm consistently performed much better than the mFLMS yielding better convergence and steady-state performance.

\bibliography{Comments_on_mFLMS}

\end{document}